\documentstyle[12pt]{amsart}

\numberwithin{equation}{section}

\renewcommand{\thesection}{\arabic{section}}
\renewcommand{\theequation}{\thesection.\arabic{equation}}

\newtheorem{theorem}[equation]{Theorem}
\newtheorem{lemma}[equation]{Lemma}
\newtheorem{corollary}[equation]{Corollary}
\newtheorem{proposition}[equation]{Proposition}

\theoremstyle{definition}
\newtheorem{definition}[equation]{Definition}

\theoremstyle{remark}

\newcounter{itemno}

\newcommand{\Pic}{\operatorname{Pic}}

\begin{document}
\title[Projective Degenerations of K3 Surfaces]
{Projective Degenerations of K3 Surfaces, \\
Gaussian Maps, and Fano Threefolds}
\author{Ciro Ciliberto}
\address{Dip. di Matematica\\Universita di Roma II\\Via Fontanile di
Carcaricola\\ 00173 Roma\\ITALY}
\email{ciliberto@@mat.utovrm.it}
\author{Angelo Lopez}
\address{Dip. di Matematica\\Universita di Pavia\\Strada Nuova 65 \\ 27100
Pavia\\ITALY}
\email{lopez@@ipvian.bitnet}
\author{Rick Miranda}
\address{Dept. of Mathematics\\Colorado State University\\Ft. Collins, CO
80523\\USA}
\email{miranda@@riemann.math.colostate.edu}
\thanks{Research supported in part by the NSF under grant DMS-9104058}

\translator{}
\dedicatory{}
\maketitle

\begin{abstract}
In this article we exhibit certain projective degenerations
of smooth $K3$ surfaces of degree $2g-2$ in $\Bbb P^g$
(whose Picard group is generated by the hyperplane class),
to a union of two rational normal scrolls,
and also to a union of planes.
As a consequence we prove that the general hyperplane section
of such $K3$ surfaces has a corank one Gaussian map,
if $g=11$ or $g\geq 13$.
We also prove that the general such hyperplane section
lies on a unique $K3$ surface, up to projectivities.
Finally we present a new approach to the classification
of prime Fano threefolds of index one,
which does not rely on the existence of a line.
\end{abstract}

{\small
\setlength{\baselineskip}{0.32cm}
\tableofcontents
}

\section*{Introduction}

Let $C$ be a stable curve of genus $g$.
There is a natural map $\phi_C$,
called the Gaussian map (or Wahl map) associated to $C$,
\[
\phi_C: \bigwedge^2 H^0(C,\omega_C) \to
H^0(C,\Omega^1_C \otimes \omega^{\otimes 2}_C)
\]
which is defined, in local coordinates, by sending $fdz\wedge gdz$ to
$(f'g-fg'){(dz)}^3$.
(Note that the target is just $H^0(C,\omega^{\otimes 3}_C)$ if $C$ is smooth.)
J. Wahl \cite{wahl} has shown that
if $C$ is the smooth hyperplane section
of a $K3$ surface of degree $2g-2$ in $\Bbb P^g$,
then the Gaussian map cannot be surjective;
see also \cite{beauville-merindol}.
More generally, F. Zak has shown that if
the canonical curve $C$ is $k$-extendable,
then the corank of the Gaussian map for $C$ is at least $k$
(see \cite{bertram-ein-lazarsfeld}).
Moreover, the corank of the Gaussian map for $C$
governs the dimension of the tangent space to the Hilbert scheme
at points representing cones over hyperplane sections
of varieties to which $C$ extends.

Thus the study of the Gaussian map for $C$
leads naturally to information concerning varieties
whose curve section is the canonical image of $C$,
namely $K3$ surfaces, Fano threefolds, etc.,
and their Hilbert schemes.
In this paper we are concerned with the following main questions:
What is the corank of the Gaussian map for the general curve $C$
which {\em is} the hyperplane section of a $K3$ surface?
What consequences can we draw for Fano threefolds
from the knowledge of this corank?

Concerning the first question we prove that,
if $g = 11$ or $g \geq 13$,
then the general canonical curve
which is the hyperplane section of a $K3$ surface,
and which generates the surface's Picard group,
has a Gaussian map with corank one (Theorem \ref{corank1_theorem}).
In genus $12$, we prove that the general such curve
has a Gaussian map with corank two (Proposition \ref{g=12corank2}).

For genus $g=10$ and $g \geq 12$,
the general curve has a surjective Gaussian map
(see \cite{ciliberto-harris-miranda}).
For  $g \leq 9$ and $g=11$ the corank of the Gaussian map
for the general curve has been studied in \cite{ciliberto-miranda1}.
In genus $10$, F. Cukierman and D. Ulmer have shown that the corank
of the Gaussian map is four for the general hyperplane section of a
$K3$ surface.

These results can be applied to the Hilbert scheme for $K3$ surfaces
and Fano threefolds.
Our main results concerning $K3$ surfaces are that,
if $g = 11$ or $g \geq 13$,
then a general hyperplane section of a $K3$ surface
(whose Picard group is generated by the hyperplane class)
lies on only one such surface, up to projectivities.
This generalizes a result of S. Mukai \cite{mukai}.
In low genera we give a new proof of a result of S. Mori and S. Mukai
(see \cite{mori-mukai}) to the effect that the general curve of genus
at most $9$ or equal to $11$ lies on a $K3$ surfaces.

Turning to the Fanos,
these techniques give a new approach
to the classification of prime Fano threefolds
which avoids completely the question of the existence of lines
and the method of double projection.
First we recover the sharp genus bound for prime Fanos,
namely, that such Fanos exist only if $g \leq 10$ or $g = 12$.
Then we prove that for these genera,
the Hilbert scheme of prime Fanos is reduced and irreducible,
and a dense open subset of it represents those prime Fanos
which have been classically exhibited by Fano (in genus up to $10$)
and by Iskovskih (in genus $12$).

Returning to the question of the corank of the Gaussian map
for hyperplane sections of $K3$ surfaces,
previous partial results in this direction
were obtained by C. Voisin \cite{voisin};
her results indicate that
the corank of the Gaussian map is at most three
for a general curve on a $K3$ surface
of high enough genus.

Our method for the proof of the corank one theorem is a degeneration technique,
similar in spirit to that of \cite{ciliberto-harris-miranda}
in that the degenerate curves are suitable graph curves.
In that work, the only difficulty is to produce graph curves
with surjective Gaussian map;
by the general theory graph curves are limits of smooth curves,
and by the semicontinuity one obtains
the surjectivity for the general smooth curve.

In this work the degeneration is quite a bit more tricky to construct.
Firstly, we must produce graph curves with corank one Gaussian map,
and then we must show that these particular stable curves are limits
of hyperplane sections of $K3$ surfaces
(with Picard group generated by the hyperplane section).
We are therefore led to degenerating not only the curves
but also the $K3$ surfaces themselves,
to suitable configurations of unions of planes,
whose hyperplane sections are then the desired graph curves.

Thus the technical part of our proof
is in producing the appropriate projective degenerations of $K3$ surfaces.
We found it convenient to make the degeneration in two steps.
We first show that the general union
of two rational normal scrolls
(each of degree $g-1$ in $\Bbb P^g$),
meeting transversally along a smooth anti-canonical elliptic curve,
can be smoothed to a $K3$ surface
which has Picard group generated by its hyperplane class.
This is described in Section \ref{defs_of_scrolls},
and uses standard deformation-theoretic arguments.
We consider this an embedded version of an abstract
Type II degeneration of $K3$ surfaces,
as described in \cite{kulikov}; see also \cite{BGOD}.
This degeneration, which to our knowledge
has not been used or systematically constructed previously,
is interesting in its own right in that it provides a purely algebraic
construction of the component of the Hilbert scheme of such $K3$ surfaces.
It is known by transcendental methods that this component is unique,
but it would be interesting to provide an algebraic proof.

Secondly we show that the union of the two scrolls
can be degenerated to a suitable union of planes
(whose dual complex is a decomposition of the sphere
with a cubic graph as its $1$-skeleton).
The particular unions of planes are described in Section \ref{union_of_planes},
and enjoy the property that the general hyperplane section
are graph curves (a union of lines) having a corank one Gaussian map.
The degeneration is constructed in Section \ref{degs_of_scrolls_to_planes}
by projective methods;
it is driven by a suitable degeneration of the elliptic double curve
and certain linear systems on this curve
which induce the scrolls.
We consider this an embedded version of an abstract
Type III degeneration of $K3$ surfaces,
as described in \cite{kulikov}.

Such degenerations are also of independent interest;
in general it is not known
whether a given configuration of planes can be smoothed.
In particular we do not know whether any configuration of planes,
whose dual complex is a decomposition of the sphere
with a cubic graph as its $1$-skeleton,
is a degeneration of a $K3$ surface.
Our work suggests that this is the case under certain extra hypotheses
(see the end of Section \ref{examples}),
but we only produce a series of examples in this article,
one for each genus.

In Section \ref{Gaussmaps} we use the degenerations constructed in the previous
sections to draw the conclusions about the Gaussian maps described above.
In this section we also make the applications
to $K3$ surfaces and Fano threefolds described above.

The authors would like to thank the University of Rome II, the C.N.R.,
and Colorado State University
for making the completion of this project possible.

\section{Deformations of the Union of Two Scrolls to Smooth $K3$ Surfaces}
\label{defs_of_scrolls}

\subsection{Smooth Elliptic Curves, $g^1_2$'s, and Rational Normal Scrolls}
\label{ENCandRNS}
Before proceeding to the main part of the section,
we need to recall some basic properties of elliptic normal curves
and rational normal scrolls.
Let $E$ be a smooth elliptic normal curve of degree $g+1$ in $\Bbb P^g$.
Let $L$ be a $g^1_2$ on $E$.
The union of the secants to the members of $L$
forms a rational normal scroll $R_L$,
which is a surface of degree $g-1$ in $\Bbb P^g$.
Moreover $E$ is an anti-canonical divisor in the scroll $R_L$.
Conversely, every rational normal scroll
containing a smooth elliptic normal curve $E$
can be constructed in this manner.
Moreover, the general scroll is obtained in this way.

In our situation we are considering the union of two scrolls
which meet along an elliptic normal curve $E$.
Given the above, we want to consider both of these scrolls
as being determined by (different) $g^1_2$'s on $E$.
The following Lemma shows that the simultaneous construction of the two
scrolls is well-behaved.

\begin{lemma}
\label{2scrolls}
Let $g$ be at least $3$, and
let $E$ be an elliptic normal curve of degree $g+1$ in $\Bbb P^g$.
Let $L_1$ and $L_2$ be two distinct $g^1_2$'s on $E$.
Let $R_1$ and $R_2$ be the corresponding scrolls.
Then $R_1$ and $R_2$ meet only along $E$,
and meet transversally along $E$.
\end{lemma}

\begin{pf}
The theorem is obvious for $g=3$, so we may assume $g \geq 4$.
First we remark that $E$ has no $4$-secant $2$-planes;
indeed, given any three points of $E$, the system of hyperplanes
through the $3$ points cuts out on $E$
a complete linear series of degree at least $2$,
and therefore has no base points.

Now suppose that there is an intersection point $p$ of $R_1$ and $R_2$
outside $E$.
There is at least one line from each of the scrolls through $p$,
and these are different since the $g^1_2$'s are different;
these two lines would span a $2$-plane, which is $4$-secant to $E$.

A similar argument shows that the intersection is transverse;
one replaces $p$ by the infinitely near point of tangency to $E$
at a supposed point of non-transversality.
\end{pf}

\subsection{The Deformation Theory for the Union of Two Scrolls}
In this section we will prove the following theorem.

\begin{theorem}
\label{smoothing_scrolls}
Fix $g \geq 3$, and let $R_1$ and $R_2$ be two smooth rational normal scrolls
each of degree $g-1$ in $\Bbb P^g$.
Assume that the two scrolls meet transversally
along a smooth elliptic normal curve $E$ of degree $g+1$,
which is anticanonical in each scroll.
Then the union $R = R_1 \cup R_2$ is a flat limit
of a family of smooth $K3$ surfaces of degree $2g-2$ in $\Bbb P^g$,
which is represented by a reduced component ${\cal H}_g$ of the Hilbert scheme
of dimension $g^2+2g+19$.
\end{theorem}

We note that the existence of an anticanonical smooth elliptic curve
in each of the two scrolls
almost implies they are general.
Indeed, in even genus they are both forced to be isomorphic to $\Bbb F_1$,
and in odd genus they are isomorphic either to
$\Bbb F_0 (\cong \Bbb P^1\times\Bbb P^1$) or to $\Bbb F_2$.

For the proof of the Theorem we need to introduce
the $T^1$ sheaf for the union $R$ of the two scrolls.
Recall that $T^1$ is defined by the sequence
\[
0 \to T_R \to T_{\Bbb P^g}|_R \to N_R \to T^1 \to 0,
\]
and for a variety with normal crossings, as $R$ is,
this sheaf $T^1$ is locally free of rank one on the singular locus.
By \cite{friedman}, we have the computation
\begin{equation}
\label{T1formula}
T^1 \cong N_{E/R_1} \otimes N_{E/R_2}
\end{equation}
where $E$ is the double curve.
In our case each sheaf $N_{E/R_i}$ is locally free rank one, of degree $8$;
therefore $T^1$ has degree $16$.

For our purposes we need the following.

\begin{lemma}
\label{F=T1}
The sheaf $\;T^1$ is isomorphic to the cokernel of the inclusion
$N_X \to N_R|_X$,
where $X$ is either $R_1$ or $R_2$.
\end{lemma}

\begin{pf}
A local computation shows that the cokernel $F$ in question
is locally free of rank one on the double curve $E$.
Consider the diagram
\[
\begin{array}{ccccccc}
T_{\Bbb P^g}|_R & \to & N_R & \to      & T^1    & \to     & 0 \\
                &     &     & \searrow &        &         &\\
\downarrow      &     &     &          & N_R|_X &         &\\
                &     &     & \nearrow &        & \searrow & \\
T_{\Bbb P^g}|_X & \to & N_X &          &        &         & 0
\end{array}
\]
where the horizontal and diagonal rows are exact.
The commutativity of the pentagon
shows that $T^1$ surjects onto $F$.
Since both are locally free rank one, this surjection must be an isomorphism.
\end{pf}

Theorem \ref{smoothing_scrolls} is, in turn,
a consequence of the following.

\begin{lemma}
\label{2scroll_lemma}
With the notation of Theorem \ref{smoothing_scrolls}, one has
\begin{itemize}
\item[(a)] $H^1(N_R) = H^2(N_R) = 0$, and
\item[(b)] $\dim H^0(N_R) = g^2 + 2g + 19$.
\end{itemize}
\end{lemma}

\begin{pf}
Because of the transversality, the normal sheaf $N_R$ is locally free.
We begin with the proof of (a).
We have the exact sequence
\[
0 \to N_R|_{R_1}(-E) \to N_R \to N_R|_{R_2} \to 0,
\]
so that it suffices to show that if we set $X = R_1$ or $R_2$,
then
\begin{equation}
\label{star}
H^j(N_R|_X) = H^j(N_R|_X(-E)) = 0 \text{ for } j = 1,2.
\end{equation}

{}From the Lemma \ref{F=T1},
to prove (\ref{star}), it therefore suffices to prove
\begin{equation}
\label{(a)}
H^j(N_X) = H^j(N_X(-E)) = 0 \text{ for } j = 1,2, \text{ and }
\end{equation}
\begin{equation}
\label{(b)}
H^j(T^1) = H^j(T^1\otimes {\cal O}_X(-E)) = 0 \text{ for } j = 1,2.
\end{equation}

\noindent
{\em Proof of \ref{(a)}:}
{}From the Euler sequence
\[
0 \to {\cal O}_X \to {{\cal O}_X(1)}^{g+1} \to T_{\Bbb P^g}|_X \to 0,
\]
it follows that $H^1(T_{\Bbb P^g}|_X) = H^2(T_{\Bbb P^g}|_X) = 0$.
Note also that $H^2(T_X) = 0$, since
\[
{H^2(T_X)}^* \cong H^0(\Omega^1_X \otimes \omega_X) = H^0(\Omega^1_X(-E))
\subseteq H^0(\Omega^1_X) = 0.
\]
Then, from the normal bundle sequence
\[
0 \to T_X \to T_{\Bbb P^g}|_X \to N_X \to 0
\]
we see that $H^j(N_X) = 0$ for $j = 1,2$.

By twisting the Euler sequence by ${\cal O}_X(-E)$,
and remarking that
\[
H^2({\cal O}_X(1)\otimes {\cal O}_X(-E))
\cong H^2(\omega_X(1))
\cong {H^0({\cal O}_X(-1))}^* = 0,
\]
we see that $H^2(T_{\Bbb P^g}(-E)|_X) = 0$;
By twisting the normal bundle sequence by ${\cal O}_X(-E)$,
we deduce that $H^2(N_X(-E))=0$.

Dualizing the normal bundle sequence and taking cohomology, we see that
\[
\begin{array}{ccccccc}
H^0(\Omega^1_X) & \to & H^1(N^*_X) & \to & H^1(\Omega^1_{\Bbb P^g}|_X) &
\stackrel{\phi}{\to} & H^1(\Omega^1_X). \\
\| & & & & \| & & \| \\
0  & & & & {\Bbb C} & & {\Bbb C}^2
\end{array}
\]
Note that the map $\phi$ is injective,
since its image is the span of the class of the hyperplane section
of $X$ in $H^{1,1}(X)$.
Therefore
\[
H^1(N_X(-E)) \cong {H^1(N^*_X)}^* = \ker \phi = 0,
\]
proving (\ref{(a)}).

\noindent
{\em Proof of \ref{(b)}:}
Of course since $T^1$ is supported on $E$,
the $H^2$'s vanish.
Moreover since $\deg(T^1) = 16$, $H^1(T^1) = 0$.
Finally $\deg(T^1\otimes {\cal O}_X(-E)) = 8$,
so that $H^1(T^1\otimes {\cal O}_X(-E))= 0$.
This finishes the proof of (\ref{(b)}),
and therefore of statement (a) of the Lemma.

We now turn to the proof of statement (b) of the Lemma.
We have that
\begin{eqnarray*}
\dim H^0(N_R) & = & \chi(N_R) \\
&=& \chi(N_R|_X) + \chi(N_R|_X(-E)) \\
&=& \chi(N_X) + \chi(T^1) + \chi(N_X(-E)) + \chi(T^1\otimes {\cal O}_X(-E)) \\
&=& g^2 +2g + 18 + \chi(N_X(-E)) \\
&=& g^2 + 2g + 18 + \chi(T_{\Bbb P^g}|_X(-E)) + \chi(T_X(-E)) \\
&=& g^2 + 2g + 18 + (g+1)\chi({\cal O}_X(1)\otimes {\cal O}_X(-E))
- \chi({\cal O}_X(-E)) + \chi(T_X(-E)).
\end{eqnarray*}
Riemann-Roch on the surface $X$ computes
$\chi({\cal O}_X(1)\otimes {\cal O}_X(-E)) = 0$
and $\chi({\cal O}_X(-E)) = 1$,
whereas $\chi(T_X(-E)) = \chi(\omega_X\otimes T_X) = \chi(\Omega^1_X) = 2$.
This finishes the proof of statement (b), and the proof of the Lemma.
\end{pf}

Note that we have as a by-product of this analysis the following.

\begin{corollary}
\label{H0NtoH0T1}
The natural map from $H^0(N_R)$ to $H^0(T^1)$ is surjective.
\end{corollary}

\begin{pf}
First note that since $H^1(N_X) = 0$ (for $X = R_1$ or $R_2$),
we have a surjection from $H^0(N_R|_X)$ to $H^0(T^1)$,
by Lemma \ref{F=T1}.
Therefore it suffices to show that $H^0(N_R)$ surjects onto $H^0(N_R|_X)$,
which is a consequence of (\ref{star}).
\end{pf}

\begin{pf*}{Proof of Theorem \ref{smoothing_scrolls}}
Because $H^1(N_R) = 0$,
$R$ is represented by a smooth point $h$ in the Hilbert scheme.
Therefore $R$ belongs to a single reduced component ${\cal H}_g$
of the Hilbert scheme of dimension $\dim H^0(N_R) = g^2 + 2g + 19$,
by Lemma \ref{2scroll_lemma}.
By Corollary \ref{H0NtoH0T1},
a general tangent vector to ${\cal H}_g$ at $h$
represents a first-order embedded deformation of $R$
which smooths the double curve.
Therefore the general point in ${\cal H}_g$
represents a smooth irreducible surface $S$.

Since $H^1({\cal O}_R) = H^1({\cal O}_R(1)) = 0$
($R$ is Cohen-Macaulay),
the same is true for $S$.
Since the hyperplane section of $R$ is a stable canonical curve of genus $g$,
the general hyperplane section of $S$ is also a smooth canonical curve.
This shows that $S$ is regular, and has trivial canonical bundle.
Therefore $S$ is a smooth $K3$ surface.
\end{pf*}

\subsection{The Picard Group of the General Smoothing}

In this subsection we want to identify the component ${\cal H}_g$
of the Hilbert scheme which we have discovered as the component
containing the union $R$ of the two scrolls.
We will prove in this section that the general surface
represented by a point in ${\cal H}_g$
has its Picard group generated by the hyperplane class.

\begin{definition}
A {\em prime} $K3$ surface of genus $g$ is a smooth $K3$ surface $S$
of degree $2g-2$ in $\Bbb P^g$
such that $\Pic(S)$ is generated by the hyperplane class.
\end{definition}

\begin{theorem}
\label{Pic=Z}
With the assumptions of Theorem \ref{smoothing_scrolls},
the general point of the component ${\cal H}_g$ of the Hilbert scheme
containing the union $R$ of the two scrolls represents a smooth
prime $K3$ surface of genus $g$.
\end{theorem}

\begin{pf}
The argument is essentially the same as that given in \cite{griffiths-harris}
and we will only outline the main steps.
We first note that we may assume
that the two $g^1_2$'s on $E$
(giving rise to the two scrolls $R_1$ and $R_2$ of $R$)
and the hyperplane class on the double curve $E$
are linearly independent in $\Pic(E)$.

Now consider a general flat family $S_t$ degenerating to $R$,
parametrized by $t$ in a disc $\Delta$;
the total space of this family in $\Bbb P^g \times \Delta$
is a singular threefold, with ordinary double points at the $16$ points
$\{p_j\}$
corresponding to the zeroes of the section of $T^1$
determined by the first-order embedded deformation.

Suppose now that we have a line bundle $W_t$ on the general fiber.
We will assume for simplicity that $W_t$ is rationally determined,
namely that there is a line bundle on the total space restricting to $W_t$
on the general fiber.  If this is not the case, one can make a further
base change and achieve this, at the cost of complicating somewhat the
singularities of the total space.  This is easily handled,
and goes exactly like the argument in \cite[Appendix C]{griffiths-harris}

Now make the total space $\cal S$ for the family smooth,
by making a small resolution of each of the $16$ double points.
This we can do in such a way that
the central fiber $S_0$ is the union of the scroll $R_1$ and
the $16$-fold blowup $\tilde{R}_2$ of $R_2$;
each of the $16$ exceptional curves $E_i$ for the double point resolutions
become $(-1)$-curves in $\tilde{R}_2$.
The two surfaces $R_1$ and $\tilde{R}_2$ still meet transversally along $E$.

We have two natural bundles on the total space $\cal S$,
namely the hyperplane bundle $H$
and the bundle $M = {\cal O}_{\cal S}(\tilde{R}_2)$.
We claim that $\Pic(S_0)$ is freely generated
by the restriction $H_0$ and $M_0$ of these two bundles.
If this is true, then we may finish the proof as follows.
Take the line bundle $W$ on the total space which restricts to $W_t$ on
the general fiber.
By twisting $W$ appropriately with multiples of $H$ and $M$,
we may assume that $W$ restricts to $H_0$ on the central fiber.
Because the general fiber is a $K3$ surface,
the dimension of $H^0(W_t)$ is numerically determined,
and so is constant over the entire family;
in particular, $h^0(W_t) = h^0(H_0) = g+1$.
Hence if $\pi$ is the map from $\cal S$ to $\Delta$,
we have that $\pi_*W$ is a trivial bundle of rank $g+1$ on $\Delta$.
Moreover $W_t$ has no base points, since the limit $H_0$ has none.
Therefore a general section of $\pi_*W$
provide us with a flat family of smooth curves in $\Bbb P^g$
of degree $2g-2$ degenerating to a
general section of $H_0$ on the central fiber,
which is a general hyperplane section of the union $R$.
This limit is a canonical curve in $\Bbb P^g$,
hence so is the general member.
Therefore the general member of the linear system for $W_t$ is a
canonical curve, which is a hyperplane section;
hence $W_t$ is the hyperplane bundle.

Alternatively, the degeneration argument shows that $W_t \cdot H_t = 2g-2$;
if $W_t$ is not isomorphic to $H_t$,
then the complete linear series cut out by $W_t$ on $H_t$
has degree $2g-2$ and dimension at least $g$.
This is forbidden by Riemann-Roch,
showing that $W_t \cong H_t$.

To finish we must prove the claim that
$\Pic(S_0)$ is freely generated by $H_0$ and $M_0$.
Now $\Pic(S_0) \cong \Pic(R_1) \times_{\Pic(E)} \Pic(\tilde{R}_2)$,
that is, a line bundle on $S_0$ is determined by line bundles on each
component agreeing along the double curve.
Let $H_i$ be the hyperplane class for $R_i$
(and $\tilde{R}_2$ if $i=2$),
and let $F_i$ be the class of the fiber in $R_i$.
Note that $F_i$ restricted to $E$ gives the corresponding $g^1_2$.

We have $\Pic(R_1)$ generated by $H_1$ and $F_1$,
and $\Pic(\tilde{R}_2)$ generated by $H_2$, $F_2$, and the $16$ exceptional
classes $E_j$.
Suppose that we have the two classes
$a_1H_1 + b_1F_1$ on $R_1$ and $ a_2 H_2 + b_2 F_2 + \sum_j c_j E_j$ on
$\tilde{R}_2$ which agree on the double curve $E$.
Then
\[
a_1 H_E + b_1 g^1_{2,1} = a_2 H_E + b_2 g^1_{2,2} + \sum_j c_j p_j
\]
must hold in $\Pic(E)$.
By a standard monodromy argument
(see \cite[Sublemma on page 36]{griffiths-harris}),
all of the coefficients $c_j$ must be equal, say to $c$.
Note that $\sum_j p_j \equiv 4H_E - (g-3)g^1_{2,1} - (g-3)g^1_{2,2}$,
by (\ref{T1formula}).
As noted above, since we are able to choose the $g^1_2$'s so
that they, along with $H$, are independent, we see by equating coefficients,
we have $a_1-a_2 = 4c$, $b_1 = (3-g)c$, and $b_2 = (g-3)c$.
This proves that $\Pic(S_0)$ depends freely on the two parameters $c$ and
$a_2$.
The bundle $H_0$ is that given by $c=0$ and $a_2 = 1$;
the bundle $M_0$ is given by $c=1$ and $a_2 = -2$.
\end{pf}

A priori there are many components to the Hilbert scheme
of prime $K3$ surfaces of genus $g$,
and what the above theorem shows is that our component ${\cal H}_g$
is one of them.
However, the transcendental theory for $K3$ surfaces
assures us that in fact there is only one component
for the Hilbert scheme of prime $K3$ surfaces of genus $g$,
which we will call the prime component.
(All other components of the Hilbert scheme of $K3$ surfaces of degree
$2g-2$ in $\Bbb P^g$ have
the Picard group of the general member
generated by a proper submultiple of the hyperplane class.)
Thus we know that our component ${\cal H}_g$ is the (unique) prime component;
but except where explicitly noted, we will not use this fact.

\section{The Union of Planes and the Graph Curves}
\label{union_of_planes}

In this section we will describe the union of planes
to which we will degenerate $K3$ surfaces in the component ${\cal H}_g$.
The general hyperplane section of the union of planes
is a union of lines, and both the union of planes and the union of lines
are conveniently described by a graph.
Indeed, these curves have been studied previously in several papers,
(see \cite{bayer-eisenbud}, \cite{ciliberto-harris-miranda}, and
\cite{ciliberto-miranda1})
and have been referred to as {\em graph curves}.

\subsection{Graph Curves and Unions of Planes}
First let us briefly recall the graph curve construction.
A graph $G$ is {\em simple} if it has no loops or multiple edges;
it is {\em trivalent} if each vertex lies on exactly three edges.
Given a simple trivalent connected graph $G$,
one forms an abstract stable curve $C_G$
by taking a smooth rational curve for each vertex of $G$,
and connecting them transversally at one point as the edges of $G$ indicate.

Since $G$ is trivalent, there is an integer $g$ such that
$G$ has $2g-2$ vertices and $3g-3$ edges.
This number $g$ is also the dimension of $H^1$ of the graph,
and is the arithmetic genus of the graph curve $C_G$.

A connected graph $G$ is {\em $n$-edge-connected}
if the removal of less than $n$ edges never disconnects $G$.
(Hence $1$-connected simply means connected;
$2$-connected means that there are no ``disconnecting'' edges of $G$.)
Recall that if $G$ is $3$-edge-connected,
then the canonical map for the graph curve $C_G$
embeds $C_G$ into $\Bbb P^{g-1}$, \cite{bayer-eisenbud};
the components of $C_G$ go to straight lines.

Now suppose that the trivalent graph $G$ is embedded in the plane,
giving not only vertices and edges, but also faces.
Using Euler's formula, one sees that there are $g+1$ faces in this
decomposition of the plane.
It is elementary to see that
if $G$ is simple, trivalent, and $2$-edge-connected,
then for each vertex, the three adjacent faces are different;
we assume this from now on.

Given the $2$-edge-connected planar graph $G$,
we form a union of $2$-planes in $\Bbb P^g$ as follows.
Take $g+1$ independent points in $\Bbb P^g$,
and give a $1$-$1$ correspondence between these points
and the faces of the decomposition.
Each vertex of $G$ lies on three faces,
which are then associated to three independent points in $\Bbb P^g$,
generating a $2$-plane.
The union of these $2$-planes over all vertices of $G$
will be denoted by $S_G$,
although it depends not only on the graph $G$ but also on the planar
embedding.

We want to have the situation that the line intersections
of the planes of $S_G$
are exactly encoded by the edges of $G$.
This is the case if and only if $G$ is $3$-edge-connected
(see \cite{ciliberto-miranda2}).
I.e., if $G$ is $3$-edge connected,
(as we will assume from now on),
then two planes of $S_G$ intersect along a line
if and only if the corresponding vertices of $G$ are joined by an edge.

Note that a general hyperplane section of $S_G$
is a canonically embedded $C_G$.
Moreover, arithmetically, the surface $S_G$ is a $K3$ surface;
namely, the Hilbert function is that of a smooth $K3$ surface
embedded in $\Bbb P^g$.
This is seen by remarking that $S_G$ is arithmetically Cohen-Macaulay
\cite[Lemma 2.12]{ciliberto-miranda1},
and has a stable canonical curve as its hyperplane section.

\subsection{Examples of Planar Graphs}
Now we will present examples of $3$-edge-connected trivalent planar graphs,
which give us via the preceding constructions,
both unions of planes and graph curves,
embedded in projective space.
We give one example for each genus $g$ at least $7$,
and we denote that graph by $G_g$;
the corresponding graph curve will be denoted by $C_g$,
and the union of planes by $S_g$.

\medskip
\noindent
{\bf The Odd Genus Case.}
Suppose that $g = 2n+1$ is odd.
Form a trivalent planar graph $G_{2n+1} = G_g$
with $2g-2 = 4n$ vertices $\{v_1,\dots,v_{4n}\}$.
There are $3g-3 = 6n$ edges of $G_{2n+1}$,
which can be described as follows.
The vertices $v_1,\dots,v_n$ are arranged in a cycle, giving $n$ edges;
similarly, the vertices $v_{n+1},\dots,v_{3n}$ lie in a cycle,
giving $2n$ additional edges;
finally, the vertices $v_{3n+1},\dots,v_{4n}$ lie in a cycle,
giving $n$ more edges.  All of these $4n$ edges are distinct.
The edge $v_i$, for $1 \leq i \leq n$, is connected to the vertex $v_{n+2i-1}$;
this adds $n$ more edges to the graph.
Similarly, the vertex $v_i$, for $3n+1 \leq i \leq 4n$,
is connected to the vertex $v_{2i-5n}$,
giving the final $n$ edges to the graph $G_{2n+1}$.

This graph $G_{2n+1}$ is planar, as can be seen from the following
construction.
Begin with $n$ points in the plane forming an $n$-gon;
these will be the vertices $\{v_1,\dots,v_n\}$.
Attach $n$ pentagons to the boundary of this $n$-gon,
each meeting the original $n$-gon in one edge,
and the neighboring pentagon on either side in one edge.
Finally attach $n$ additional pentagons,
each meeting two of the previous set of pentagons,
and again meeting the neighboring pentagon on either side in one edge.
We refer the reader to Figure 1.

\medskip
\noindent
{\bf The Even Genus Case.}
The graph in the case that $g$ is even
is obtained from the graph $G_{g-1}$ constructed above
by adding two vertices and one edge.
One of the two additional vertices is added on the edge joining
vertex $1$ to vertex $n$;
the other is added on the opposite side of the central $n$-gon,
on the edge joining vertex $[(n+1)/2]$ to vertex $[(n+3)/2]$.
The added edge joins these two additional vertices.

This graph is also planar; the construction is similar to the odd genus case.
One begins instead with two polygons,
one with $[(n+4)/2]$ sides, and one with $[(n+5)/2]$ sides,
joined at one side.
The union of the two form an $(n+2)$-gon with one interior edge.
Attach $n$ pentagons to the boundary of the $(n+2)$-gon,
with exactly two of the pentagons having one of their edges covering
two of the edges of the $(n+2)$-gon,
which the interior edge meets.
Finally attach $n$ additional pentagons,
as in the odd genus case.
We refer the reader to Figure 1.
Note that in this construction $g = 2n+2$.
Also note that the numbering of the vertices is slightly different
than in the odd genus case.

In the next section we will describe a deformation of the surface $S_g$
to a union of two rational normal scrolls.
Each of the scrolls degenerates to exactly half of the $2g-2$ planes.
We will now describe the two sets of $g-1$ planes
to which the scrolls individually degenerate.

Recall that the vertices of the graph $G_g$
correspond to the $2$-planes of $S_g$.
Therefore we may describe the desired subset of planes
by indicating which of the vertices correspond to that subset.
We will refer to one subset of $g-1$ vertices as the $\cal A$ subset,
and the other as the $\cal B$ subset;
we will refer to the corresponding union of planes by $S_A$ and $S_B$
respectively.
We again have two cases, depending on the parity of the genus $g$.

\medskip
\noindent
{\bf The Odd Genus Decomposition.}  Recall here that $g = 2n+1$,
so that each subset $\cal A$ and $\cal B$ consist of $2n$ vertices.
We set:
\[
\cal A = \{ v_1, v_{n+1},\dots,v_{3n-2}, v_{4n-1} \}, \text{ and }
\]
\[
\cal B = \{ v_2,\dots,v_n,v_{3n-1},v_{3n},v_{4n},v_{3n+1},\dots,v_{4n-2} \}.
\]

\medskip
\noindent
{\bf The Even Genus Decomposition.} Recall here that $g = 2n+2$,
so that each subset $\cal A$ and $\cal B$ consist of $2n+1$ vertices.
We set:
\[
\cal A = \{ v_{n+1}, v_1,v_{n+3},v_{n+4},\dots,v_{3n},v_{4n+1} \},\text{ and }
\]
\[
\cal B = \{ v_2,v_3,\dots,v_{[(n+1)/2]},v_{n+2},v_{[(n+3)/2]},\dots, v_n,
v_{3n+1},v_{3n+2},v_{4n+2},v_{3n+3},v_{3n+4},\dots,v_{4n}  \}.
\]

The reader may consult Figure 2 for a representation of these decompositions.

The order chosen above for the subsets $\cal A$ and $\cal B$
may seem a bit strange;
however they are such that in that order,
the vertices for both subsets form a path in the graph $G_g$.
Since the vertices form a path,
the corresponding sets of planes $S_A$ and $S_B$ are arranged also in a path,
with each plane meeting two others along lines,
except the first and last.
This configuration is again Cohen-Macaulay,
and the hyperplane section is a degenerate form of a rational normal curve,
namely is a union of $g-1$ straight lines in $\Bbb P^{g-1}$.
Hence each of the configurations $S_A$ and $S_B$ of planes
is numerically a rational normal scroll,
in the sense that it has the same Hilbert function.

Of course, the union $S_A \cup S_B$ is the total configuration $S_g$.
The two configurations intersect along a cycle of lines;
the number of lines is $g+1$.
We have drawn on Figure 2
a dotted path separating the two subsets of vertices $\cal A$ and $\cal B$;
this path cuts the graph in exactly those edges which correspond
to the lines in which the configurations $S_A$ and $S_B$ meet.
Since the path meets each face exactly once,
this set of $g+1$ lines forms a cycle
through the coordinate points of $\Bbb P^g$
which were chosen to correspond to the faces of the planar decomposition.
In particular, this intersection is a cycle of lines
which is maximally embedded, and forms a degenerate elliptic normal curve
in $\Bbb P^g$ of degree $g+1$.

It is useful to have in mind the way in which the separate configurations
$S_A$ and $S_B$ of planes are put together.
This we show in Figure 3, for the odd genus case.
Note that the numbering of the planes is of course
the numbering of the corresponding vertices.

For later use, an alternative description of these unions of planes
will be needed.
The intersection of the two chains of planes,
as noted above, forms a cycle of $g+1$ lines.
Number these (modulo $g+1$) in order around the cycle,
and denote them therefore by $\ell_1,\ell_2,\dots,\ell_{g+1}$.
Also number the intersection points of these lines in a standard way,
so that they are denoted by $p_j$, where $p_j = \ell_j \cap \ell_{j+1}$.

Now note that each plane in the configurations
is the span of either two intersecting lines $\ell_j$ and $\ell_{j+1}$
(this is the case for the planes at the ends of the chains)
or is the span of a line $\ell_j$ and a point $p_k$.
Thus to describe the configurations one can simply give
these spanning sets. This we now do for the configuration $S_A$
in the odd genus case, assuming that the lines are numbered
as indicated in Figure 3.

\begin{center}
\begin{tabular}{||l|l||}
\multicolumn{2}{c}{Table One} \\ \hline
Plane & Span of  \\ \hline
$1$ & $\ell_1$, $\ell_{2n+2}$ \\ \hline
$n+1$ & $p_1$, $\ell_{2n+1}$ \\ \hline
$n+2i$, $1 \leq i \leq n-2$ & $p_i$, $\ell_{2n+1-i}$ \\ \hline
$n+2i+1$, $1 \leq i \leq n-2$ & $p_{2n-i}$, $\ell_{i+1}$ \\ \hline
$3n-2$ & $p_{n+2}$, $\ell_n$ \\ \hline
$4n-1$ & $\ell_{n+1}$, $\ell_{n+2}$ \\ \hline
\end{tabular}
\end{center}

\section{Degenerations of the Union of Two Scrolls to the Union of Planes}
\label{degs_of_scrolls_to_planes}

In this section we will describe a deformation
of the union of planes $S_g$,
which were constructed in the previous section,
to a union of
two rational normal scrolls.
As mentioned above, each of the configurations $S_A$ and $S_B$
will be smoothed to a scroll; care must be taken that this is possible
in such a way that the intersection curve also deforms correctly.

Let us denote the two scrolls by $R_A$ and $R_B$;
$R_A$ degenerates to $S_A$ and $R_B$ to $S_B$.
The intersection of the scrolls will be a smooth elliptic normal curve
which is an anti-canonical divisor in each scroll;
this intersection curve will degenerate to the cycle of lines
which is the intersection of the $S_A$ and $S_B$ configuration of planes.
Our first step is to reduce the analysis of the degeneration
to these intersection curves.

The construction of scrolls as the secants
to a $g^1_2$ on an elliptic normal curve
as described in Section \ref{ENCandRNS}
reduces our degeneration construction to the construction
of an appropriate degeneration of the elliptic normal curve,
along with the degeneration of the two $g^1_2$'s.

\subsection{Cycle Degenerations of Elliptic Normal Curves}
\label{cycle_degenerations}
Let us next discuss the degeneration of an elliptic normal curve
which we will use.
For every integer $k \geq 2$,
there is an abstract degeneration of smooth elliptic curves
to a cycle of $k$ rational curves, namely the degeneration
denoted by $I_k$ in Kodaira's original construction (see \cite{kodaira}).
The total space of this degeneration is a smooth surface $T_k$,
mapping to the disc $\Delta$,
and the general fiber of the map is an elliptic curve.
The central fiber is a cycle of $k$ $\Bbb P^1$'s,
each of self-intersection $-2$.
This degeneration is unique in the sense that
any two such degenerations are isomorphic in a neighborhood of the
central fiber.
Moreover, one can arrange the isomorphism
to identify one component of the central fiber of one
with any component of the central fiber of the other,
and also to make the identifications of the two cycles in
any of the two orderings for the cycle.

Suppose that $k \geq g+1$.
Choose any $g+1$ of the $k$ components
of the central fiber of the degeneration $T_k$.
In addition, choose $g+1$ sections of the map to $\Delta$,
each meeting one of the chosen $g+1$ components.
These sections give a divisor of degree $g+1$ in each fiber
of the map to $\Delta$,
and the corresponding linear system embeds the fibers into $\Bbb P^g$
as projectively normal curves,
giving a map, defined over $\Delta$, from $T_k$ to $\Bbb P^g \times \Delta$.
The smooth fibers go to elliptic normal curves of degree $g+1$ in $\Bbb P^g$.
The central fiber goes to a cycle of $g+1$ lines which span $\Bbb P^g$.
In particular, each of the $g+1$ components which are met by the $g+1$ sections
map to one of the lines in the cycle;
the $k-g-1$ components not meeting any of the sections
are contracted (to rational double points of type $A$).

\subsection{The Degeneration of the $g^1_2$'s}
For the degeneration of the $g^1_2$'s,
we need to view a $g^1_2$ on an elliptic curve
not as a linear system
but as a $2$-$1$ branched cover of $\Bbb P^1$.
With this point of view,
a degeneration of a $g^1_2$
becomes a double covering of an appropriate degeneration of $\Bbb P^1$'s.
This idea is essentially the notion of ``admissible coverings'',
in a very special case.

Firstly we present the relevant construction of such degenerations
of $\Bbb P^1$ and the double covering.
Start with a trivial family $\Bbb P^1 \times \Delta \to \Delta$.
Blow up in the central fiber $m-1$ times,
creating an abstract degeneration $P_m$ of $\Bbb P^1$
to a chain of $m$ $\Bbb P^1$'s;
the ends of the chain are both $(-1)$-curves,
and the interior members of the chain are all $(-2)$-curves.
In order to have a double covering of $P_m$ whose
general fiber is elliptic, the branch locus must meet the general fiber
in $4$ points.  Choose a smooth $4$-section $B$ which meets
the central fiber only in the two end components,
and meets each of these twice.
Let $D_m(B)$ be the double covering of the surface $P_m$;
it is a smooth surface, mapping to $\Delta$,
and is a degeneration of elliptic curves.

There are two cases for the behavior of the branch curve $B$
at the two ends of the chain: $B$ could either meet the end component
at two points transversally, or $B$ could be tangent to the end component.
Moreover either of these cases can happen at either end of the chain.
This gives rise to three cases:
\begin{itemize}
\item[(i):] {\em The branch curve $B$ meets both end components
transversally at two points each.}\\
In this case the double covering $D_m(B)$ is
a degeneration of elliptic curves isomorphic to $T_{2m-2}$.
Each of the interior $m-2$ components of $P_m$
splits into disjoint curves,
and the two end components are doubly covered by a single component.
\item[(ii):] {\em The branch curve $B$ meets one component
transversally at two points,
and is tangent to the other component at one point.}\\
In this case the double covering $D_m(B)$ is
a degeneration of elliptic curves isomorphic to $T_{2m-1}$.
Each of the interior $m-2$ components of $P_m$ again split,
and the end component which $B$ meets transversally is also doubly covered
by a single component.
However the other end component to which $B$ is tangent also splits,
into two components meeting transversally
at the point lying over the point of tangency.
\item[(iii):] {\em The branch curve $B$ is tangent to both end components.}\\
In this case the double covering $D_m(B)$ is
a degeneration of elliptic curves isomorphic to $T_{2m}$.
Each of the interior $m-2$ components of $P_m$
again split into disjoint curves,
while the two end components to which $B$ is tangent split into two curves
which meet transversally as above.
\end{itemize}

The double covering $D_m(B) \to P_m$
induces an involution on each fiber of $D_m(B)$ over $\Delta$.
On the central fiber, this involution becomes
a correspondence between the components.
In case (i), the components which cover the end components
carry a self-involution;
the other components split into two chains (of length $m-2$),
and the involution on the fiber makes a correspondence
between the two chains of components.
In case (ii), only the single component of $D_m(B)$
covering the component of $P_m$
which meets the branch curve $B$ in two points
carries a self-involution;
the remaining $2m-2$ components are paired.
Finally, in case (iii),
the $2m$ components are all paired, in two chains of $m$ components,
starting and ending with the components covering the end components of $P_m$.

These three situations lead us to define the following notion.

\begin{definition}
\label{allowable_double_correspondence}
Let $\cal C$ be a cycle of smooth rational curves of length $k$.
An {\em allowable double correspondence} on $\cal C$ is one of the following
data:
\begin{itemize}
\item[(i):] If $k$ is even, with $k = 2m-2$,
one gives a component of $\cal C$, its ``opposite'' component on the cycle,
and a correspondence between the remaining two chains of components
which associate components with the same distance
from the two chosen components.
In particular, if one numbers the components $C_i$ around the cycle
so that the chosen components are $C_1$ and $C_m$,
then in the two chains, the components $C_i$ and $C_{2m-i}$ are associated.
\item[(ii):] If $k$ is odd, with $k = 2m-1$,
one gives a component of $\cal C$,
and a correspondence between the remaining components
which associates components with the same distance
from the chosen component.
In particular, if one numbers the components $C_i$ around the cycle
so that the chosen component is $C_1$,
then the components $C_i$ and $C_{2m-i+1}$ are associated.
\item[(iii):]If $k$ is even, with $k = 2m$,
one gives a partition of the components of $\cal C$
into two chains of length $m$,
and a correspondence between these chains,
which associate components with the same distance
from one of the two vertices of the cycle in which the two chains intersect.
In particular, if one numbers the components $C_i$ around the cycle so that
the first component on one of the two chains is $C_1$,
then the components $C_i$ and $C_{2m-i+1}$ are associated.
\end{itemize}
\end{definition}

We remark that most of this data is actually determined by a single choice.
In case (i), if one chooses one of the components of the cycle $\cal C$,
then the opposite component and the correspondence between the remaining
components are all determined.  Similarly in case (ii),
it suffices to choose a single component, and in case (iii),
a single vertex.

It is obvious that for the degenerations of the $g^1_2$'s
which we constructed above, each case in the construction
(that is, each case for the branch locus $B$),
gives the corresponding allowable double correspondence on the central fiber.
Moreover, we have a converse to this construction.

\begin{proposition}
\label{adcprop}
Let $X \to \Delta$ be an $I_k$ degeneration of elliptic curves,
that is, $X$ is a smooth surface mapping to the disc $\Delta$ with
general fiber an elliptic curve and with central fiber a cycle of $k$
$\Bbb P^1$'s, each having self-intersection $-2$.
Suppose that on the central fiber there is given
an allowable double correspondence.
Then $X$ is isomorphic (in a neighborhood of the central fiber)
to $D_m(B)$ for the appropriate $m$ and $B$,
in such a way that the allowable double correspondences coincide.
\end{proposition}

\begin{pf}
By the uniqueness of the $I_k$ degenerations,
we know that $X \to \Delta$ is isomorphic to $D_m(B) \to \Delta$
for the appropriate $m$ and $B$.
Moreover, as we mentioned above, given isomorphic $I_k$ degenerations,
one can arrange the isomorphism to make any two components correspond.
By adjusting the components which are associated via the isomorphism,
we may clearly make the two allowable double correspondences to coincide.
\end{pf}

\subsection{Admissible Projective Embeddings of Degenerations of $g^1_2$'s}

In this subsection we will combine the ideas in the previous subsections
and explain how, in certain cases, with a degeneration of $g^1_2$'s
and a suitable projective embedding, one can see the degeneration of the
induced rational normal scrolls to a union of planes.

The complete situation is the following.
Suppose we are given an $I_k$ degeneration of elliptic curves
$X \to \Delta$, abstractly isomorphic to the $T_k$ degeneration,
with $k \geq g+1$.
In addition, suppose that exactly $g+1$ of the components
of the central fiber have been chosen,
which induces a cycle degeneration of $X$ into $\Bbb P^g \times \Delta$
as described at the end of subsection \ref{cycle_degenerations}:
the $g+1$ chosen components survive in the limit to form a cycle of lines.
Finally we assume that an allowable double correspondence has been given
on the central fiber of $X$.

By Proposition \ref{adcprop},
we have a degeneration of $g^1_2$'s on $X$,
which gives us on the general fiber $X_t$ of the map to $\Delta$
a $g^1_2$ determining  a rational normal scroll
$R_t$ in $\Bbb P^g \times \{t\}$.
This is a flat family of surfaces over the punctured disc,
and we want to understand the limit of these scrolls in  $\Bbb P^g$.

While this question is interesting in general,
in our application we need only to understand some particular cases
for the relationship between the choice of $g+1$ surviving components
and the allowable double correspondence.
These we describe next.

\begin{definition}
\label{compat}
A choice of $g+1$ surviving components is {\em compatible} with
an allowable double correspondence on the central fiber of a $T_k$
degeneration of elliptic curves if and only if
the following conditions are satisfied:
\begin{itemize}
\item[(a)] In case (i) or (ii), no component which carries the self-involution
can be one of the $g+1$ surviving components, that is, these end components
must be contracted in the map to projective space.
\item[(b)] In the two chains of components
which are paired in the allowable double correspondence,
the first component to survive (on either side) is paired with a surviving
component on the other chain.
\item[(c)] Except for these two pairs of surviving components,
no other surviving component is paired with a surviving component.
\end{itemize}
\end{definition}

This definition is exactly that which makes the analysis of the limit of the
scrolls rather transparent.

\begin{proposition}
\label{planes}
Suppose that a set of $g+1$ surviving components is chosen
on the degeneration $T_k \to \Delta$,
which is compatible with an allowable double correspondence.
Then the flat limit of the rational normal scrolls $R_t$
is the union of $g-1$ $2$-planes,
which are spanned by the following.
\begin{itemize}
\item[($\alpha$)] For each of the two pairs of surviving components
as in (b) above, these components are mapped to lines which meet at one point;
therefore they span a $2$-plane, which is part of the limit.
\item[($\beta$)] For each of the other $g-3$ surviving components,
which map to a line, the paired component is contracted to a point not on that
line; this point and line span a $2$-plane, which is part of the limit.
\end{itemize}
\end{proposition}

\begin{pf}
Since the scroll is determined
as the union of secants to the general $g^1_2$,
the limit of the scrolls will contain the secants to the limit $g^1_2$.
Recall that in our situation the degeneration of the $g^1_2$'s
is obtained by the double covering of the central fiber,
as in the example $D_m(B) \to P_m$.
By definition, the secants in the limit are made between points of the cycle
on $D_m(B)$ which map to the same point in the double covering.
Therefore, in the limit, the secants are made between points of the paired
components in the allowable double correspondence.

When neither of a pair survives, we have no $2$-dimensional contribution
to the limit.
When both survive, as in case ($\alpha$), we clearly obtain the plane
spanned by the components.
When exactly one of the pair survives, as in case ($\beta$),
we again have the plane spanned by the surviving component and by the
point to which the contracted component is mapped.

The union of these obvious $g-1$ $2$-planes contained in the limit
is a reducible surface composed of a chain of $2$-planes,
and is arithmetically a rational normal scroll, that is,
it has the same Hilbert function as a rational normal scroll.
This proves that this union of planes is exactly the flat limit of the
scrolls $R_t$.
\end{pf}

\subsection{Examples}
\label{examples}

In this subsection we will exhibit the degeneration of the union of
two scrolls in $\Bbb P^g$ to the union of $2g-2$ planes,
in the configuration of the planar graphs given in the examples
of Section \ref{union_of_planes}.
Our method will be to give a $T_k$ degeneration of elliptic curves,
together with a choice of $g+1$ surviving components,
and two different allowable double correspondences on the central fiber,
which are both compatible with the choice of surviving components.
By Proposition \ref{planes}, this induces two different degenerations
of scrolls through the general elliptic curve,
both of which degenerate to planes in a well-controlled manner.
We will show that the limit is exactly the union of the planes
$S_A \cup S_B$ which we desire.
Since the $g^1_2$'s in the limit are different,
the $g^1_2$'s on the general fiber must be different,
and hence the two scrolls $R^{(1)}_t$ and $R^{(2)}_t$
in the general fiber will meet transversally along the general elliptic curve
by Lemma \ref{2scrolls}.

The union of the planes $S_A \cup S_B$
is indeed the flat limit of the two scrolls $R^{(1)}_t$ and $R^{(2)}_t$.
In fact, both $S_A \cup S_B$ and $R^{(1)}_t \cup R^{(2)}_t$
have the same Hilbert function, namely that of a smooth $K3$ surface
of degree $2g-2$ in $\Bbb P^g$.

\begin{theorem}
\label{2scroll_limit}
For each $g \geq 7$, the union of $2g-2$ $2$-planes $S_g$
described in Section \ref{union_of_planes}
is a flat limit of a union of two rational normal scrolls,
meeting transversally along a smooth elliptic curve,
which is anti-canonical in each scroll.
\end{theorem}

As noted above, to prove the Theorem it is sufficient to exhibit
the appropriate allowable double correspondences,
which are compatible with a choice of $g+1$ surviving components
on some $T_k$ degeneration of elliptic curves.
The construction naturally breaks up into the odd and even genus cases,
and we will treat these separately.

\medskip
\noindent
{\bf The Odd Genus Case.}
Recall that in this case we set $g = 2n+1$;
now also set $k = 4g-8 = 8n-4$.
This is the length of the cycle before the embedding into projective
space.
Number the components around the cycle from $1$ to $k=8n-4$,
consecutively.
As noted above we must give the $g+1 = 2n$ surviving components,
as well as two different allowable double correspondences,
which will give us the $\cal A$ and $\cal B$ configurations of planes
in the limit.
Let $\cal G$ be the set of indices of the surviving components.
Then in this case:
\[
\cal G = \{n-2,n+3,n+5,n+7,\dots,3n+1,5n-4,5n+1,5n+3,5n+5,\dots,7n-1\}.
\]

Now we give the two allowable double correspondences,
which we also denote by $\cal A$ and $\cal B$.
The numbering has been chosen so that $\cal A$ is the pairing
of component $j$ with component $8n-3-j$:
\[
\cal A: j \leftrightarrow 8n-3-j;
\]
this is the numbering used to describe case (iii) of the allowable
double correspondence definition, in fact.
The other allowable double correspondence $\cal B$
is also of type (iii),
and is the pairing of component $j$ with component $2n+1-j \mod 8n-4$:
\[
\cal B: j \leftrightarrow 2n+1-j \;  \mod 8n-4.
\]

Let us check that $\cal G$ is compatible with $\cal A$;
we will leave the analogous proof for $\cal B$ to the reader.
Firstly, both $\cal A$ and $\cal B$ are of type (iii),
so condition (a) of Definition \ref{compat} is automatic.
Secondly, the first surviving components on the two chains of $\cal A$
at one end are indexed $n-2$ and $7n-1$, and these are paired;
similarly, at the other end, the first surviving components are
indexed $3n+1$ and $5n-4$, and these are paired.
Thus condition (b) is satisfied.
Finally, the indices of the interior components
which survive on one side are all of the
form $n + \text{ odd number }$,
and these are paired with components whose indices are of the form
$7n - \text{ even number }$; no such component survives in $\cal G$.
Thus condition (c) of Definition \ref{compat} is also satisfied,
and we have shown that $\cal G$ is compatible with $\cal A$.

If we renumber the surviving components from $1$ to $g+1 = 2n+2$,
simultaneously numbering the vertices in the standard way,
then the planes of the limit of the $\cal A$ $g^1_2$,
according to Proposition \ref{planes}
clearly form the configuration $S_A$
described in Section \ref{union_of_planes}
(in particular, in Table One) for this genus.
The same is true for the $\cal B$ $g^1_2$:
it has as the limit the configuration $S_B$.
This completes the proof of Theorem \ref{2scroll_limit}
in the odd genus case.

\medskip
\noindent
{\bf The Even Genus Case.}
Recall that in this case we set $g = 2n+2$;
now also set $k = 8n + 6[(n-1)/2] - 11$.
This is the length of the cycle before the embedding into projective
space.
Number the components consecutively around the cycle from $0$ to $k-1$.
We must again give the $g+1 = 2n+3$ surviving components,
as well as two different allowable double correspondences,
which will give us the $\cal A$ and $\cal B$ configurations of planes
in the limit.
Let $\cal G$ be the set of indices of the surviving components.
Then in this case:
\begin{eqnarray*}
\cal G &=& \{ [(n-2)/2],[(3n-7)/2],
[(3n+1)/2],[(3n+5)/2],[(3n+9)/2],\dots,[(7n-7)/2], \\
& & 3[(3n-3)/2],2n+3[(3n-3)/2]-5, \\
& & 5[(3n-3)/2],5[(3n-3)/2]+2,5[(3n-3)/2]+4,\dots, \\
& & 8n+4[(n-1)/2]-[(n-2)/2]-11,8n+6[(n-1)/2]-[(n-2)/2]-11 \}
\end{eqnarray*}

Now we give the two allowable double correspondences,
which we also denote by $\cal A$ and $\cal B$.
The numbering has been chosen so that $\cal A$ is the pairing
of component $j$ with component $k-j$:
\[
\cal A: j \leftrightarrow k-j \; \mod k;
\]
note that $0$ is self-paired, and is the only such index,
so that $\cal A$ is of type (ii) in the allowable
double correspondence definition.
The other allowable double correspondence $\cal B$
is also of type (ii),
and is the pairing of component $j$ with component $2[(3n-3)/2]-j \mod k$:
\[
\cal B: j \leftrightarrow 2[(3n-3)/2]-j \;  \mod k.
\]

We leave to the reader in this even genus case
to check that $\cal G$ is compatible with $\cal A$ and $\cal B$,
and that the limit configuration of planes
is exactly the union $S_A \cup S_B$ for this genus.

This completes the proof of Theorem \ref{2scroll_limit}.

As a consequence of Theorem \ref{2scroll_limit},
Theorem \ref{smoothing_scrolls}, and Theorem \ref{Pic=Z},
by combining the deformations of the planes to the two scrolls,
and then the deformation of the two scrolls to the $K3$ surfaces,
we have the following.

\begin{corollary}
\label{main_theorem}
For each $g \geq 7$,
the union of planes $S_g$ described in Section \ref{union_of_planes}
is a flat limit of smooth prime $K3$ surfaces of genus $g$,
and in particular is represented by a point of ${\cal H}_g$.
Therefore, the graph curve $C_g$,
which is the general hyperplane section of the union of planes $S_g$,
is a flat limit of hyperplane sections of such smooth prime $K3$ surfaces
of genus $g$.
\end{corollary}

These examples suggest that, in general, if the union of the $2g-2$ planes
determined by a planar graph
can be decomposed into two sets of $g-1$ planes,
each of which form a chain of planes,
then the union can be smoothed to a union of two scrolls.

\subsection{Remarks for Low Values of the Genus}
\label{low_genus_remarks}

Actually, for low genus, it is not difficult to degenerate smooth $K3$
surfaces to unions of planes, whose dual graph is planar.
For $g=3$, for example, one has the quartic surface degenerating to
the union of the four coordinate planes in $\Bbb P^3$.
Also for $g=4,5$, the general smooth $K3$ surface is a complete intersection,
and the degeneration to the union of planes is straightforward.

Up to genus $10$ we have also the following construction.
There are smooth $K3$ surfaces in this range,
obtained by intersecting a cone over a Del Pezzo surface
with a general quadric.
If one allows the Del Pezzo surface to degenerate to a cone over
a cycle of lines, and the quadric to degenerate to a union of hyperplanes,
then the smooth surface degenerates to a union of planes,
whose dual graph is a ``prism'',
that is, consists of two cycles of vertices joined in series.
(In genus $10$ the smooth $K3$ surfaces are not prime,
but they are prime for genus up to $9$.)

Since the analysis is trivial in genus less than $6$,
we now present another alternate degeneration
in genus $7$ and $8$, taken from the planar graphs
introduced in \cite{ciliberto-miranda1}.
These graphs give graph curves with the minimal possible corank
of the Gaussian map, which is neither the case for the prism graphs,
nor for the graphs $G_7$ and $G_8$,
and therefore it is interesting to note that
the corresponding unions of planes
are limits of $K$ surfaces.

For completeness we set $\tilde{G}_6$, $\tilde{G}_7$, and $\tilde{G}_8$
to be the planar graphs depicted in \cite[page 429]{ciliberto-miranda1}.
The graph $\tilde{G}_6$ is a prism over a pentagon,
and can therefore be smoothed as described above.

For the union of planes $\tilde{G}_7$,
we return to the construction via the degeneration of the elliptic curve
and the two $g^1_2$'s.  Set $k=18$, and consider the $T_k$
degeneration of elliptic curves, numbered in order around the cycle from
$1$ to $18$.  The set $\cal G$ of surviving components is
\[
\cal G = \{1,5,7,10,11,13,15,16\}.
\]
The allowable double correspondence giving the $\cal A$ set of planes
is of type (iii), and pairs component $j$ with component $17-j\mod 18$.
The allowable double correspondence giving the $\cal A$ set of planes
is of type (i), with the components $3$ and $12$ being self-paired,
and otherwise pairs component $j$ with component $6-j\mod 18$.

For the union of planes $\tilde{G}_8$,
we set $k=23$, and consider the $T_k$
degeneration of elliptic curves,
numbered in order around the cycle from $1$ to $23$.
The set $\cal G$ of surviving components is
\[
\cal G = \{1,6,8,9,12,14,16,19,20\}.
\]
The allowable double correspondence giving the $\cal A$ set of planes
is of type (ii), and pairs component $j$ with component $21-j\mod 23$;
component $22$ is self-paired.
The allowable double correspondence giving the $\cal B$ set of planes
is also of type (ii),
with component $15$ being self-paired,
and otherwise pairs component $j$ with component $7-j\mod 23$.

The reader may check that these are compatible allowable double
correspondences, and give as the limit of the two scrolls,
the union of planes described by the graphs $\tilde{G}_j$,
for $j=7,8$.  Indeed, the reader may reconstruct these graphs now
from the data given above, by working out the union of planes.

Finally we remark that the graph given in \cite{ciliberto-miranda1}
for genus $9$ is the same as the $G_9$ graph of this paper.

\section{Applications to Gaussian Maps, $K3$ Surfaces, and Fano Threefolds}
\label{Gaussmaps}

\subsection{The Corank One Theorem for $g=11$ and $g \geq 13$}
\label{corank_one_section}

Let $C$ be a stable curve.
There is a natural map
\[
\phi:\bigwedge^2H^0(\omega_C) \to H^0(\Omega^1_C\otimes\omega^{\otimes 2}_C)
\]
defined for a smooth curve in local coordinates by
$\phi(fdz\wedge gdz) = (f'g-fg'){(dz)}^3$.

Note that the corank of the Gaussian map is semi-continuous in moduli.
Indeed, the question is local
on the compactified moduli space $\overline{{\cal M}}_g$,
and so one may work over the Kuranishi family for the stable curve.
In this case, by constancy of dimensions and standard base change theorems,
the domain and range of the Gaussian map
fit together to form a vector bundle over the base,
and the map is a map of bundles.  The semi-continuity follows.

Wahl \cite{wahl} has shown that a canonical curve
which is the hyperplane section of a $K3$ surface
cannot have a surjective Gaussian map (see also \cite{beauville-merindol}).
Using the degeneration techniques developed in this article,
we can make this theorem more precise.
Let ${\cal H}_g$ be the component of the Hilbert scheme of $K3$ surfaces
of degree $2g-2$ in $\Bbb P^g$ which degenerate to the union of two scrolls,
as described in Theorem \ref{smoothing_scrolls}.

\begin{theorem}
\label{corank1_theorem}
Suppose that $g=11$ or $g\geq 13$.
Let $S$ be a $K3$ surface of degree $2g-2$ in $\Bbb P^g$
represented by a general point in ${\cal H_g}$.
Then the general hyperplane section $C$ of $S$
has a Gaussian map with corank exactly one.
\end{theorem}

\begin{pf}
First we remark that since the corank of the Gaussian map is semi-continuous,
and it cannot be surjective by Wahl's result,
it suffices to prove the theorem for a single stable curve $C$
which is the limit of general hyperplane sections of such $K3$ surfaces.
For this purpose we use the graph curves $C_g$,
which are the hyperplane sections of the union of planes $S_g$
introduced in Section \ref{union_of_planes}.
These are limits of such $K3$ surfaces,
as we have shown in Corollary \ref{main_theorem}.
The Gaussian maps for such graph curves
have been the subject of several papers by the authors.
In particular, in \cite{miranda},
it is shown that each of these graph curves,
in genus $11$ or genus at least $13$,
have a corank one Gaussian map.
This is also proved in \cite{ciliberto-franchetta1},
where a more systematic study is made and more general results are obtained.
In any case this computation completes the proof of the Theorem.
\end{pf}

L. Ein has communicated to us an alternate approach to the corank one theorem.
His approach would first prove the surjectivity
of the relevant Gaussian map defined on the $K3$ surface itself;
this would imply that every smooth hyperplane section
would have a corank one Gaussian map,
not simply the general one.
His technique requires a decomposition of the hyperplane bundle
which seems difficult to realize in low genera.
However Ein's approach gives sharp results
for hyperplane sections of $K3$ surfaces which are re-embedded
via Veronese maps.
We intend to treat this in a forthcoming paper.

\subsection{The Cases of Genus Twelve and Ten}
\label{g=12case}

In genus $12$, the general curve has a surjective Gaussian map
(see \cite{ciliberto-harris-miranda}).
For the graph curve $C_{12}$ constructed in this article,
and which is a limit of $K3$ sections by Corollary \ref{main_theorem},
the Gaussian map has corank two,
by the computations presented in \cite{ciliberto-franchetta2}.
Therefore by semi-continuity,
the general $K3$ section in genus $12$ has corank either one or two.
In fact it is two.

We require a lemma.

\begin{lemma}
\label{h0N(-k)}
Let $C$ be a general hyperplane section of a general prime $K3$ surface of
genus $g$, with $g \geq 6$.
If $k \geq 3$, then $H^0(N_C(-k)) = 0$.
If $g \geq 7$, then $H^0(N_C(-2)) = 0$;
if $g = 6$, then $\dim H^0(N_C(-2)) \leq 1$.
When $k=1$, we have $\dim H^0(N_C(-1)) = g + \gamma_C$,
where $\gamma_C$ is the corank of the Gaussian map for $C$.
\end{lemma}

\begin{pf}
For $k \geq 3$, $H^0(N_C(-k)) = 0$ since the ideal of $C$ is generated
by quadrics. For $k=2$, a computation on the graph curves
introduced in this paper
(namely those associated to the graphs $\tilde{G}_g$ for $g=7,8$
and the graphs $G_g$ for $g \geq 9$)
shows that $H^0(N_C(-2))=0$ for $g \geq 7$.
The computation is exactly that of \cite[Theorem 3.1]{ciliberto-miranda1};
in fact, this space is zero for any graph curve associated to a graph
which is not a prism.  The genus $6$ case was also computed in
\cite[Theorem 3.1]{ciliberto-miranda1}, and the above statement
follows from semi-continuity.
The final statement concerning $H^0(N_C(-1))$ is
\cite[Proposition 1.2(a)]{ciliberto-miranda2}.
\end{pf}

\begin{proposition}
\label{g=12corank2}
The general hyperplane section of a $K3$ surface
of degree $22$ in $\Bbb P^{12}$
represented by a general point in ${\cal H}_{12}$
has corank two Gaussian map.
\end{proposition}

\begin{pf}
We must only show that the general such hyperplane section
cannot have corank one.
By \cite{mukai}, the generic prime $K3$ surface of genus $12$
is the hyperplane section of a smooth Fano $3$-fold.
(This is true in lower genera also; we discuss this in general below.)
Therefore the general $K3$ hyperplane section $C$ is $2$-extendable,
i.e., is a codimension two linear section of a smooth variety.
If the corank of the Gaussian map for $C$ is one, then by Lemma \ref{h0N(-k)}
we have $H^0(N_C(-2)) = 0$ and $\dim H^1(N_C(-1)) = g + 1$.
These are the hypothesis of a theorem of Zak/L'vovsky
(see \cite{lvovsky} or \cite{bertram-ein-lazarsfeld}),
which concludes that $C$ cannot be $2$-extendable.
Since this is a contradiction,
we see that the corank of the Gaussian map must be two.
\end{pf}

By \cite{cukierman-ulmer}, the general hyperplane section
of a prime $K3$ surface
of genus $10$ has a corank $4$ Gaussian map.
Using the techniques of this paper,
it might be possible to prove,
by exhibiting an example of such a curve,
that the corank of the Gaussian map
for the general prime $K3$ section of genus $10$ was at most four.
We have been unable to construct an appropriate planar graph curve
of genus $10$ with corank $4$ Gaussian map.
We leave this as an open problem: to find such a graph curve.

The cases of lower genera,
in particular $g \leq 9$ and $g=11$,
are discussed rather fully in \cite{ciliberto-miranda1}.

Note that the above statements rely on the transcendental theory
for $K3$ surfaces, which is used to insure that our component ${\cal H}_g$
is the prime component.

\subsection{The Family of $K3$ Sections}
\label{families}

Let ${\cal H}_g$ be the component of the Hilbert scheme of $K3$ surfaces
of degree $2g-2$ in $\Bbb P^g$ which degenerate to the union of two scrolls,
as described in Theorem \ref{smoothing_scrolls}.
In this and the next section
we will invoke the transcendental theory
and refer to this component ${\cal H}_g$
as the Hilbert scheme for prime $K3$ surfaces of genus $g$.

Let ${\cal F}_g$ be the ``flag Hilbert scheme''
parametrizing pairs $(S,C)$,
where $S$ is represented by a point of ${\cal H}_g$
and $C$ is a stable hyperplane section of $S$.
Finally let ${\cal C}_g$ be the Hilbert scheme of degenerate
stable canonical curves of genus $g$ in $\Bbb P^g$,
which of course all live in some hyperplane.
Note that
\begin{eqnarray*}
\dim {\cal H}_g &=& g^2 + 2g + 19, \\
\dim {\cal C}_g &=& g^2 + 4g - 4, \text{ and }\\
\dim {\cal F}_g &=& g^2 + 3g + 19.
\end{eqnarray*}
We have natural maps $p:{\cal F}_g \to {\cal C}_g$,
and $q:{\cal F}_g \to {\cal H}_g$.
In this subsection we will study some properties of the map $p$.
We will fix the genus $g$ to be at least $6$;
the lower genera are those of complete intersections,
and all of the analysis is trivial.

First let us compute the dimension of the general fiber of $p$.
By standard deformation theory (see \cite{kleppe} or \cite{kleppe2}),
the tangent space to the fiber at a point $(S,C)$ is isomorphic to
$H^0(N_S \otimes I_{C/S})$, where $N_S$ is the normal bundle to $S$
and $I_{C/S}$ is the ideal sheaf of $C$ in $S$.

We apply this remark to the pair $(X,C)$,
where $C$ is a general point in the image of $p$,
(i.e., $C$ is a general prime $K3$ surface section),
and $X$ is a general cone over $C$.
Since any projectively Cohen-Macaulay surface
flatly degenerates to a cone over its hyperplane section,
$X$ is indeed represented by a point in ${\cal H}_g$
and the pair $(X,C)$ by a point $(x,c)$ in ${\cal F}_g$.

For the cone $X$, $H^0(N_X) \cong \oplus_{k\geq 0} H^0(N_C(-k))$;
hence the tangent space to the fiber of $p$ at the point $(x,c)$
is isomorphic to $\oplus_{k\geq 1} H^0(N_C(-k))$.

Hence, using Lemma \ref{h0N(-k)},
the tangent space to the fiber of $p$ at $(x,c)$
has dimension $g+\gamma_C$ for $g \geq 7$,
and is at most $17$ for genus $6$.

Now the discussion breaks into several cases.
First assume that $6 \leq g \leq 9$, or $g = 11$.
In this case we have an upper bound for the corank of the Gaussian
map, by the graph curve computations made in
\cite{ciliberto-miranda1};
this gives, by the above Lemma, and semi-continuity,
an upper bound for the dimension of the general fiber of $p$.
The answer (see \cite[Theorem 3.3]{ciliberto-miranda1} is that
the fiber dimension is at most $23-g$.
Hence in this range, we have that $p$ is surjective;
in particular, the general fiber dimension is exactly $23-g$,
and the inequality above in the genus $6$ case is an equality.
The surjectivity of $p$ is also proved in \cite{mori-mukai},
in a different manner.

We notice that, even ignoring the transcendental theory,
one has that ${\cal H}_g$ is the only component
of the Hilbert scheme of $K3$ surfaces which surjects onto $\cal C_g$,
as proved in \cite{ciliberto-miranda1}.

Next assume that $g \geq 13$.
For this range, we know that the corank of the Gaussian map
is one for the general curve $C$ in the image of $p$, by
Theorem \ref{corank1_theorem}.
Hence the tangent space to the fiber of $p$ at the cone point $(x,c)$
is $g+1$, and hence by semi-continuity,
the tangent space to the fiber of $p$ at a general point is at most $g+1$.
On the other hand, the fiber clearly has dimension at least $g+1$,
by using projective transformations fixing the hyperplane section $C$.
Therefore the general fiber has dimension $g+1$.

Finally let us take up the cases of genus $12$ and $10$.
In genus $12$, by the above Lemma and Proposition \ref{g=12corank2},
the dimension of the general fiber of $p$
is at most $14$.  In fact it is exactly $14$,
since the generic prime $K3$ surface of genus $12$
is the hyperplane section of a smooth Fano threefold,
as described in \cite{mukai}.
In genus $10$ the same remarks apply:
using the results of \cite{cukierman-ulmer}
the dimension of the general fiber of $p$ is again $14$.

These computations are sufficient to show the following.

\begin{proposition}
\label{fiber_irred}
Fix $g \geq 3$.  Then the general fiber of the map $p$ is irreducible.
\end{proposition}

\begin{pf}
Suppose that there are two components in the fiber.
Then the point $(x,c)$,
where $c$ is a general point in the image of $p$
and $x$ represents a cone over the curve represented by $c$,
belongs to both components,
and hence is a singular point of the fiber.
However, we have shown above that the tangent space to the cone point
is exactly the dimension of the fiber.  Thus the fiber must be irreducible.
\end{pf}

The above Proposition was proved for low genera in
\cite[Theorem 5.3]{ciliberto-miranda1}.

Let $F_g$ be the moduli space for prime $K3$ surfaces of genus $g$.
Let $I_g$ be the moduli space of pairs $(S,C)$,
where $(S,{\cal O}_S(C)) \in F_g$,
and $C$ is a stable curve on $S$.
We have a natural projection $\psi:I_g \to F_g$,
and a natural map $\pi: I_g \to {\cal M}_g$,
where ${\cal M}_g$ is the moduli space for curves of genus $g$.
We note that $\dim F_g = 19$, and $\dim I_g = g+19$.

\begin{theorem}
\label{moduli_map}
For $3 \leq g \leq 9$ and for $g=11$,
the map $\pi$ is dominant, and the general fiber is irreducible.
For $g=10$, the codimension of the image of $\pi$ is one,
and for $g=12$, the codimension of the image of $\pi$ is two.
For $g \geq 13$, the map $\pi$ is birational onto its image.
Moreover, for $g = 11$ and $g \geq 13$,
the general canonical curve
which is the hyperplane section of a prime $K3$ surface
lies on a unique one, up to projective transformations.
\end{theorem}

The proof is an easy consequence of the previous Proposition,
since the fibers of the map $\pi$ at the moduli space level
are obtained by taking the fibers of the map $p$ at the Hilbert scheme level
and dividing by the action of the projective group.

\begin{corollary}
\label{codim1conelocus}
Let ${\cal H}_g$ be the Hilbert scheme of prime $K3$ surfaces of genus $g$.
Let ${\cal X}_g$ be the locus in ${\cal H}_g$ representing cones over
hyperplane sections of the surfaces in ${\cal H}_g$.
Then the codimension of ${\cal X}_g$ in ${\cal H}_g$ is
$\gamma_C + \sum_{k\geq 2}h^0(N_C(-k))$
where $C$ is a general curve in the image of $p$.
In particular, ${\cal X}_g$ has codimension one in ${\cal H}_g$
if and only if $g=11$ or $g \geq 13$.
\end{corollary}

\begin{pf}
Let $\delta$ be the dimension of the image of $p$.
By the above discussion,
we have that the general fiber of $p$ has dimension
$g+\gamma_C + \sum_{k\geq 2}h^0(N_C(-k))$
where $C$ is a general curve in the image of $p$.
Since the map $q$ is surjective with $g$-dimensional fibers,
we see that ${\cal H}_g$ has dimension
$\delta + \gamma_C + \sum_{k\geq 2}h^0(N_C(-k))$.
On the other hand, it is elementary that ${\cal X}_g$ has dimension $\delta$.
This proves the general statement,
and the final statement follows from Theorem \ref{corank1_theorem}.
\end{pf}

We remark that if one considers a general point of ${\cal X}_g$,
representing a cone $X$ with vertex $v$,
and if we are in the corank one case of $g=11$ or $g\geq 13$,
then the only embedded deformations of $X$ are either to smooth $K3$ surfaces
(represented by points of ${\cal H}_g$)
or to other cones.
Therefore this singularity cannot be ``partially smoothed'':
any deformation either is a total smoothing,
or is topologically trivial.

In conclusion, let us remark that
the results here lead naturally to the following questions.
Suppose that $g \geq 13$,
and consider the moduli space ${\cal M}_g$ for curves of genus $g$.
Let ${\cal N}_g$ be the locus representing curves
with non-surjective Gaussian maps.
Let ${\cal K}_g$ be the locus representing hyperplane sections
of smooth prime $K3$ surfaces of genus $g$.
Is ${\cal K}_g$ a component of ${\cal N}_g$?
What are the components of ${\cal N}_g$ when it is not irreducible?
The first question has a positive answer in genus $10$,
by \cite{cukierman-ulmer},
and in this case ${\cal N}_g$ is irreducible.
Finally, a related question is:
does every stable curve with corank one Gaussian map
lie on a numerical $K3$ surface?
We know of no counterexamples to this; the other examples
of curves with corank one Gaussian map are planar graph curves,
or hypersurface sections of cones over canonical curves,
all of which lie on numerical $K3$'s.

\subsection{Fano Threefolds}
\label{fanos}
We note that Theorem \ref{moduli_map}
extends the result of Mori and Mukai
(see \cite{mori-mukai} and \cite{mukai})
which only proves that $\pi$ is generically finite if $g \geq 13$.
Our more precise statement can be applied to the existence and irreducibility
of families of Fano threefolds as we now explain.

A {\em prime} Fano threefold of genus $g$
is a smooth projective threefold
of degree $2g-2$ in $\Bbb P^{g+1}$
whose hyperplane class is anticanonical,
and generates the Picard group.
Note that by standard Noether-Lefschetz-type arguments,
the general hyperplane section
of a prime Fano threefold of genus $g$
is a prime $K3$ surface of genus $g$
(see e.g. \cite{moishezon}).
Therefore the general curve section of a prime Fano threefold
is a canonical curve which is at least $2$-extendable;
hence the corank of the Gaussian map for these curve sections
must be at least two, by Zak's Theorem (see \cite{bertram-ein-lazarsfeld}).

Let ${\cal V}_g$ be the Hilbert scheme of prime Fano threefolds
of genus $g$.
Next we remark that the general hyperplane section $S$ of a prime Fano $V$
represented by a general point in {\em any} component of ${\cal V}_g$
describes a general prime $K3$ surface.
In other words, the second projection
from any component of the Hilbert scheme of pairs $(V,S)$
to ${\cal H}_g$ is surjective.
This follows by considering the two maps
$H^1(T_V) \to H^1(T_V|_S)$ and $H^1(T_S) \to H^1(T_V|_S)$.
The second map is surjective:
the next term in the cohomology sequence is $H^1(N_{S/V})$,
and since $N_{S/V} \cong {\cal O}_S(1)$, this $H^1$ vanishes.

The first map has corank at most one,
since the next term in the cohomology sequence is $H^2(T_V(-1))$,
which is Serre dual to $H^1(\Omega_V^1)$,
and therefore has dimension one since $V$ is prime.
On the other hand this first map cannot be surjective;
if it were, every deformation of $S$ to first order
would come from deforming $V$,
and since there are non-algebraic deformations of $S$,
this is not possible.

Hence first map has corank exactly one and the second map is surjective,
which, by standard deformation theory arguments,
proves the above statement:
a general hyperplane section of a general Fano
(in any component of ${\cal V}_g$)
is a general prime $K3$ surface.

The following theorem now follows from the above remarks
and Theorem \ref{corank1_theorem}.

\begin{theorem}
\label{nofanos}
The Hilbert scheme ${\cal V}_g$ is empty for $g=11$ and $g \geq 13$.
\end{theorem}

The above theorem was first proved by V.A. Iskovskih in \cite{iskovskih2},
and reproved in the above spirit by S. Mukai in \cite{mukai}.

Now let us turn to the question of irreducibility
of the Hilbert scheme ${\cal V}_g$.

We assume from now on that $g\geq 6$,
since the lower genera are complete interesections and the theory is trivial.
Now any prime Fano of genus $g$ can be flatly degenerated to a cone over
its general hyperplane section.
Moreover we know that ${\cal H}_g$ is the only component of the Hilbert scheme
whose general point represents a prime $K3$ surface of genus $g$.
Hence, by the same argument as was used
in the proof of Proposition \ref{fiber_irred},
it will suffice to show that the cone
over a $K3$ surface represented by a general point in ${\cal H}_g$
is represented by a smooth point of ${\cal V}_g$.

Let $X$ be such a cone over a $K3$ surface $S$ with curve section $C$.
The tangent space to the Hilbert scheme containing $X$
has dimension
$h^0(N_X) = \sum_{k\geq 0}h^0(N_S(-k)) = h^0(N_S)+h^0(N_S(-1))+h^0(N_S(-2))$
since $S$ is cut out by quadrics.
Now $h^0(N_S)$ is the dimension of ${\cal H}_g$, which is $g^2+2g+19$.
The space $H^0(N_S(-1))$ is the tangent space to the fiber of the map
$p:{\cal F}_g \to {\cal C}_g$ introduced in Section \ref{families},
and by the results there we deduce that this space has dimension
$h^0(N_C(-1))+h^0(N_C(-2))$.
Finally $h^0(N_S(-2))$ is,
by semi-continuity upon deforming to a cone over $C$,
at most $h^0(N_C(-2))$.
Therefore, by Lemma \ref{h0N(-k)},
we have the following formula:
\begin{equation}
\label{H0NX}
h^0(N_X) \leq \begin{cases}
g^2 + 3g + 19 + \gamma_C & \text{ if }g\geq 7,\text{ and }\\
g^2 + 3g + 19 + \gamma_C+ 2 & \text{ if }g=6.
\end{cases}
\end{equation}.

By the results of Fano and Iskovskih,
for $6 \leq g \leq 10$ and $g=12$,
(see \cite{iskovskih2} and \cite{mukai}),
the Hilbert scheme ${\cal V}_g$ is non-empty;
indeed, there are examples for these genera
which provide families with the number of moduli and parameters
shown in Table Two;
also shown there is the corank $\gamma_C$ of the Gaussian map
for the general curve section $C$.

\begin{center}
\begin{tabular}{||c|c|c|c||}
\multicolumn{4}{c}{Table Two} \\ \hline
genus $g$ & number of moduli & number of parameters & $\gamma_C$ \\ \hline
 6 & 22 &  85 & 10\\ \hline
 7 & 18 &  98 &  9 \\ \hline
 8 & 15 & 114 &  7 \\ \hline
 9 & 12 & 132 &  5 \\ \hline
10 & 10 & 153 &  4 \\ \hline
12 &  6 & 201 &  2\\ \hline
\end{tabular}
\end{center}

The number of moduli for each of the examples
are computed easily using projective arguments.
The number of parameters for the family in $\Bbb P^{g+1}$
which these examples fill out
is obtained simply by adding the number of projective transformations
of $\Bbb P^{g+1}$, which is of course $g^2 + 4g +3$.
The final column giving the value of the corank of the Gaussian map
is taken from \cite[Theorem 3.2]{ciliberto-miranda1}
and Section \ref{g=12case}.

Now comparing the number of parameters for these families
with the tangent space to the Hilbert scheme at the general cone,
we see that the upper bound for the tangent space dimension
given in (\ref{H0NX})
coincides for each of these genera
with the number of parameters.
Therefore we can conclude the following.

\begin{theorem}
\label{Vgirreducible}
For $6 \leq g \leq 10$ and $g=12$,
the point of ${\cal V}_g$ represented by
a cone over a general prime $K3$ surface $S$ in ${\cal H}_g$
is a smooth point of ${\cal V}_g$.
Moreover, ${\cal V}_g$ is irreducible,
and the examples of Fano and Iskovskih
fill out all of ${\cal V}_g$.
Finally, the fiber of the second projection from pairs $(V,S)$
to ${\cal H}_g$ is irreducible;
in genus $12$, this implies that the general prime $K3$ surface
lies on a unique prime Fano, up to projective transformations.
\end{theorem}

Note that the above theorem
can be viewed as a classification theorem for prime Fano threefolds.
This approach to the theory of Fano threefolds
avoids completely the necessity of proving the existence of lines
and the method of double projection.
It is clear that using this approach the entire theory and classification
of Fano threefolds can be re-formulated.
For the non-prime Fanos,
a similar analysis can be made;
this involves a careful study of the corank of the Gaussian map
and the cohomology of the normal bundle
for curves on general non-prime $K3$ surfaces.
This we will return to in a later paper,
as mentioned at the end of Section \ref{corank_one_section}.


\begin{thebibliography}{XXXNNN}
\bibitem[B-E]{bayer-eisenbud} D. Bayer and D. Eisenbud:
``Graph Curves'', preprint.
\bibitem[B-E-L]{bertram-ein-lazarsfeld} A. Bertram, L. Ein, and R. Lazarsfeld:
``Surjectivity of Gaussian Maps for Line Bundles of Large Degree on Curves''.
preprint
\bibitem[B-M]{beauville-merindol} A. Beauville and J.-Y. Merindol:
``Sectiones Hyperplanes des Surfaces $K3$'',
Duke Mathematical Journal, Vol. 55 (1987), 873-878.
\bibitem[BGOD]{BGOD}{\em The Birational Geometry of Degenerations},
edited by R. Friedman and D. Morrison.  Progress in Mathematics Vol. 29,
Birkh\"auser 1983.
\bibitem[C-F]{ciliberto-franchetta1} C. Ciliberto and A. Franchetta:
``Curve poligonali, grafi e applicazione di Gauss'',
Rendiconti di Matematica, vol. 12 (1992), 165-196.
\bibitem[CF2]{ciliberto-franchetta2}C. Ciliberto and A. Franchetta:
``L'applicazione di Gauss per grafi trivalenti'', preprint.
\bibitem[C-H-M]{ciliberto-harris-miranda}
C. Ciliberto, J. Harris, and R. Miranda:
``On the surjectivity of the Wahl map'',
Duke Mathematical Journal, vol. 57, no. 3, (1988), 829-858.
\bibitem[C-M1]{ciliberto-miranda1} C. Ciliberto and R. Miranda:
``On the Gaussian map for canonical curves of low genus'',
Duke Mathematical Journal, vol. 61, No.2 (1990), 417-443.
\bibitem[C-M2]{ciliberto-miranda2} C. Ciliberto and R. Miranda:
``Graph Curves, Colorings, and Matroids'',
to appear in the Proceedings of the Ravello Conference
on Zero-Dimensional Schemes, Ravello, Italy, June 1992.
\bibitem[C-U]{cukierman-ulmer}F. Cukierman and D. Ulmer:
``Curves of Genus Ten on $K3$ Surfaces''.
To appear in Compositio Math.
\bibitem[F]{friedman}R. Friedman:
``Global Smoothings of Varieties with Normal Crossings'',
Annals of Mathematics, Vol. 118 (1983), 75-114.
\bibitem[G-H]{griffiths-harris} P. Griffiths and J. Harris:
``On the Noether-Lefschetz Theorem and
Some Remarks on Codimension-two Cycles'',
Math. Annalen, vol. 271 (1985), 31-51.
\bibitem[I]{iskovskih2}V.A. Iskovskih:
``Fano Threefolds II'', Math. USSR Izvestija, Vol. 12 (1978), 469-506.
\bibitem[Kl1]{kleppe} J. Kleppe:
``The Hilbert Flag Scheme, its properties, and its connection
with the Hilbert Scheme. Applications to curves in three-space''.
Thesis, University of Oslo, Preprint No. 5 (1981).
\bibitem[Kl2]{kleppe2} J. Kleppe:
``Non-reduced components of the Hilbert scheme of smooth space curves'',
in: Space Curve, Rocca di papa, 1985,
ed. by F. Ghione, C. Peskine, and E. Sernesi.
Springer Lecture Notes in Mathematics, No. 1266 (1987), 181-207.
\bibitem[Ko]{kodaira} K. Kodaira:
``On compact  analytic  surfaces  II,  III".
Annals of Math. 77, (1963), 563-626; 78 (1963), 1-40.
\bibitem[Ku]{kulikov} V. Kulikov:
``Degenerations of $K3$ surfaces and Enriques Surfaces'',
Math. USSR Izvestija, vol. 11 (1977), 957-989.
\bibitem[L]{lvovsky} S. M. L'vovsky:
``On the extension of varieties defined by quadratic equations'',
Math. USSR Sbornik, Vol. 63 (1989), 305-317.
\bibitem[Mi]{miranda} R. Miranda:
``The Gaussian map for certain planar graph curves'',
In: Algebraic Geometry: Sundance 1988. B. Harbourne and R. Speiser, editors.
Contemporary Mathematics 116 (1991), 115-124.
\bibitem[Mo]{moishezon} B. Moishezon:
``On algebraic cohomology classes on algebraic varieties'',
Izv. Akad. Nauk USSR, Ser. Math 31 (1967), No. 2, 225-268.
\bibitem[M-M]{mori-mukai} S. Mori and S. Mukai:
``The uniruledness of the moduli space of curves of genus $11$''.
In: Algebraic Geometry, Proceedings, Tokyo/Kyoto (1982),
Springer LNM No. 1016 (1983), 334-353.
\bibitem[Mu]{mukai} S. Mukai: ``Fano $3$-folds''.
In: {\em Complex Projective Geometry}, edited by G. Ellingsrud,
C Peskine, G. Sacchiero, and S.A. Stromme,
London Mathematical Society Lecture Note Series No. 179 (1992), 255-263.
\bibitem[V]{voisin}C. Voisin:
``Sur l'application de Wahl des courbes
satisfaisant la condition de Brill-Noether-Petri'', preprint.
\bibitem[Wa]{wahl}J. Wahl:
``The Jacobian Algebra of a graded Gorenstein singularity".
Duke Math. Journal, vol. 55 (1987), 843-871.
\end{thebibliography}
\end{document}